%% file: main.tex
\documentclass[fleqn,5p]{elsarticle}

\journal{Applied Energy}

\usepackage{multicol}
\usepackage{amsmath}
\usepackage{amssymb}
\usepackage{mathtools}
\usepackage{cuted}
\usepackage{graphicx}
\usepackage{flafter}
\usepackage{placeins}
\usepackage{IEEEtrantools}
\usepackage{textcomp}
\usepackage[hidelinks,hyperindex,breaklinks]{hyperref}
\usepackage[hyphenbreaks]{breakurl}
\usepackage{float}
  \usepackage{tikz}
  \usepackage{pgfplots,pgfplotstable}
  \usetikzlibrary{pgfplots.groupplots}
  \pgfplotsset{compat=1.5}
\usepackage{empheq} 
\usepackage{color}
\newcommand{\norm}[1]{\left\lVert#1\right\rVert} 
\hypersetup{urlcolor=cyan}

\usepackage{flushend} 

\begin{document}

\setlength{\tabcolsep}{6pt}
\setlength{\mathindent}{0pt}

\begin{frontmatter}

\title{Pricing inertia and Frequency Response with diverse dynamics\\in a Mixed-Integer Second-Order Cone Programming formulation}
\tnotetext[mytitlenote]{
This paper was submitted for review on 14/09/2019. A portion of this work has been supported by the Engineering and Physical Sciences Research Council under grants EP/R045518/1 and EP/L019469/1.}

\author[imperial]{L. Badesa}
\ead{luis.badesa@imperial.ac.uk}
\author[imperial]{F. Teng}
\author[imperial]{G. Strbac}
\address[imperial]{Imperial College London, Department of Electrical and Electronic Engineering, London, SW7 2AZ, UK}

\begin{abstract}
Low levels of system inertia in power grids with significant penetration of non-synchronous Renewable Energy Sources (RES) have increased the risk of frequency instability. The provision of a certain type of ancillary services such as inertia and Frequency Response (FR) is needed at all times, to maintain system frequency within secure limits if the loss of a large power infeed were to occur. In this paper we propose a frequency-secured optimisation framework for the procurement of inertia and FR with diverse dynamics, which enables to apply a marginal-pricing scheme for these services. This pricing scheme, deduced from a Mixed-Integer Second-Order Cone Program (MISOCP) formulation that represents frequency-security constraints, allows for the first time to appropriately value multi-speed FR.  
\end{abstract}


\end{frontmatter}

\section{Nomenclature}


\subsection{Acronyms}
\begin{tabular}{p{12mm} l}
ED &Economic Dispatch.\\
EFR &Enhanced Frequency Response.\\
FR &Frequency Response.\\
KKT &Karush Kuhn Tucker.\\
MISOCP &Mixed-Integer Second-Order Cone Program.\\
MSG &Minimum Stable Generation.\\
SOC &Second-Order Cone.\\
PFR &Primary Frequency Response.\\
RES &Renewable Energy Sources.\\
RoCoF &Rate-of-Change-of-Frequency.\\
q-s-s &Quasi-steady-state frequency.\\
UC &Unit Commitment.
\end{tabular}

\subsection{Indices}
\begin{tabular}{p{12mm} l}
$i$ &All-purpose index.\\
$g$ & Index for thermal generators.
\end{tabular}

\subsection{Sets}
\begin{tabular}{p{12mm} l}
$\mathcal{G}$ &Set of thermal generators.
\end{tabular}

\subsection{Functions}
\begin{tabular}{p{12mm} l}
$\Delta f(t)$ &Post-fault frequency deviation (Hz).\\
$\textrm{FR}(t)$ &Aggregated system FR (MW).\\
$t_\textrm{nadir}$ &Time when the frequency nadir occurs (s).
\end{tabular}

\subsection{Decision Variables\\\normalfont{ (continuous unless otherwise indicated)}}
\begin{tabular}{p{12mm} l}
$\lambda_1,\lambda_2,\mu$ &Dual variables for the SOC nadir\\ 
                          &constraint.\\
$\lambda_\textrm{R}$ & Dual variable for the RoCoF constraint.\\
$\lambda_\textrm{qss}$ & Dual variable for the q-s-s constraint.\\
$P_g$ &Power produced by generator $g$ (MW).\\ 
$P_{\textrm{L}}$ &Largest power infeed (MW).\\
$P^\textrm{curt}_{\textrm{R}}$ &RES power curtailed (MW).\\
$R_i$ &FR provision from service $i$ (MW).\\
$R_g$ &FR provision from generator $g$ (MW).\\
$y_g$ &Binary variable, commitment decision\\
      &for generator $g$.
\end{tabular}

\subsection{Linear Expressions\\\normalfont{ (linear combinations of decision variables)}}
\begin{tabular}{p{12mm} l}
$H$ &System inertia (MW$\cdot \textrm{s}$).
\end{tabular}

\subsection{Parameters}
\begin{tabular}{p{15mm} l}
$\Delta f_{\textrm{max}}$ & Maximum admissible frequency \\
      &deviation at nadir (Hz).\\
$\textrm{c}^{\textrm{m}}_g$ &Marginal cost of generator $g$ (\pounds/MWh).\\
$\textrm{c}^{\textrm{nl}}_g$ &No-load cost of generator $g$ (\pounds/h).\\
$f_0$ &Nominal frequency of the grid (Hz).\\
$\textrm{H}_g$ &Inertia constant of generator $g$ (s).\\
$\textrm{H}_\textrm{L}$ & Inertia constant of the generator\\
                        &producing power $P_{\textrm{L}}$ (s).\\
$\textrm{P}_{\textrm{D}}$ & Total demand (MW).\\
$\textrm{P}_{\textrm{g}}^{\textrm{max}}$ & Maximum power output of generator $g$\\
                                         &(MW).\\
$\textrm{P}_{\textrm{g}}^{\textrm{msg}}$ & Minimum stable generation of\\
                                         &generator $g$ (MW).
\end{tabular}
\begin{tabular}{p{15mm} l}
$\textrm{P}_{\textrm{L}}^{\textrm{max}}$ & Upper bound for the largest power\\
                                         &infeed (MW).\\
$\textrm{P}_{\textrm{R}}$ & Power produced from RES (MW).\\
$\textrm{R}_{\textrm{g}}^{\textrm{max}}$ & FR capacity of generator $g$ (MW).\\
$\textrm{RoCoF}_{\textrm{max}}$ & Maximum admissible RoCoF (Hz/s).\\ 
$\textrm{T}_i$ &Delivery time of FR service $i$ (s).\\
$\textrm{T}_{\textrm{del},i}$ &Delay in provision of FR service $i$ (s).
\end{tabular}

\section{Introduction} \label{Intro}

The increasing penetration of Renewable Energy Sources (RES) in power grids has brought attention to ``frequency services", which are ancillary services for the compliance with frequency regulation. In order to avoid the unexpected tripping of Rate-of-Change-of-Frequency (RoCoF) relays or Under-Frequency Load Shedding, system frequency must be contained within a narrow band around the nominal value of 50Hz, in the event of any credible contingency. It is precisely these frequency services that allow to keep frequency within secure margins: system inertia, the different types of Frequency Response (FR), and a reduced largest loss through part-loading large generating units or cross-country interconnectors \cite{LuisPESGM2018}.

While the cost of operating a traditional power system, dominated by thermal generators, was mainly driven by the cost of burning fuel for energy production, modern RES-dominated power systems show a different behaviour: the overall cost of operating the grid is reduced, since there is no fuel-cost associated with RES, but the fraction of total cost due to providing frequency services has greatly increased, and could reach more than 10\% of overall cost in some cases \cite{FeiISGT2017}. The higher need for frequency services is driven by the low system inertia levels caused by the non-synchronous RES, which do not contribute to system inertia as they are decoupled from the grid through power electronics converters.

The need for finding methods to efficiently provide frequency services is clear by the ample literature available on the topic. The impact on FR requirements due to penetration of non-synchronous wind and solar generation was analysed in \cite{AppEn4}. The authors in \cite{AppEnFreq1} studied secondary FR coordination from industrial loads and thermal generators, while \cite{AppEnFreq2} focused on the coordination of electric vehicles and thermal units to provide response. Studies on the provision of FR from novel devices are also widely available, such as \cite{AppEnFreq3} for battery storage and \cite{AppEnFeiEV} for electric vehicles.

A day-ahead market for FR services has recently been envisioned by National Grid in Great Britain \cite{NationalGridRoadmap}. Although the scope of this market is to provide FR independently of the energy market, cost efficiency would be significantly enhanced by using simultaneous clearing of energy and FR \cite{GoranBook}, as well as by introducing novel products such as extra inertia and a reduced largest-loss. Note that the simultaneous clearing of energy and reserve has already been advocated by several works such as \cite{AppEnReserve}. Rewarding inertia provision has been demonstrated by several studies \cite{ElaI,LuisISGT2017} to provide significant economic benefits to the system, and a synchronous-inertial-response product has been recently introduced by EirGrid in Ireland \cite{EirGridInertiaProduct}. Optimally scheduling the largest possible loss has also been shown to achieve both economic and emission savings \cite{LuisPESGM2018,OMalleyDeload}.

In this context, some studies have focused on analysing market design for frequency services. The two-part paper \cite{ElaI,ElaII} considered a single Primary Frequency Response (PFR) service, while the work in \cite{ERCOT_EFR} and \cite{FeiCambridge} included also Enhanced Frequency Response (EFR). Previous works mainly apply regression techniques to obtain linear frequency-security constraints based on dynamic-simulation data \cite{ElaI,ElaII} or system operational data \cite{ERCOT_EFR}. From the linear constraints, they developed a marginal pricing mechanism for FR services. On the other hand, the authors in \cite{FeiCambridge} linearised an originally non-linear constraint in order to clear a market considering both the EFR and PFR services. Regarding a market for inertia, reference \cite{ElaI} used a similar approach as for its proposed market for PFR, by fitting heuristic linear constraints from simulation data and applying a marginal-pricing scheme to the linear constraint. Reference \cite{ETHthesisInertia} only focused on pricing inertia from the RoCoF constraint, for which an analytical linear formulation can be obtained. 

In this paper, we consider a purely analytical methodology to obtain marginal prices for the different frequency services, namely inertia, different types of FR services with distinct dynamics, and a reduced largest loss. The pricing mechanism is deduced from the constraints that guarantee frequency security, i.e. guarantee that RoCoF, frequency nadir and frequency quasi-steady-state (q-s-s) are within acceptable limits. Our formulation can also consider any combination of activation delays for the different FR services, and the impact of this delay on the marginal price is explicitly reflected. The impact of FR-delays on frequency stability has been studied in several recent works. The authors in \cite{QitengDelays} analysed the impact of FR-delays in post-fault frequency dynamics, by conducting a sensitivity analysis using dynamic simulations. Reference \cite{DelaysEFR} considered the delays in the FR dynamics of battery storage, while \cite{DelaysDSR} focused on the FR dynamics from demand-side response. The work in \cite{DelayStability} developed a delay heuristic that estimates the post-fault frequency excursion. However, no study has yet analysed the economic impact of delays on FR markets.

While most of existing energy markets rely on Linear programming, for example through DC Optimal Power Flow, the frequency-security constraints deduced in this paper are formulated as a Mixed-Integer Second-Order Cone Program (MISOCP). Second-Order Cones (SOCs) are the highest class of conic optimisation problems that can be efficiently solved to global optimality in a mixed-integer formulation using commercial solvers \cite{BoydConvex}. Therefore, an MISOCP formulation is suitable for being implemented in a scheduling algorithm that needs binary variables to model the on/off state of the generators, such as the Unit Commitment.

The contributions of this work are: 
\begin{enumerate}
	\item We propose a marginal pricing scheme for frequency services (multi-speed FR, inertia and a reduced largest-loss), based on duality theory applied to the MISOCP formulation that we deduce for guaranteeing post-fault frequency security. This pricing scheme allows for the first time to optimally value the different dynamics of FR providers, including the possibility to consider certain FR services with activation delays, as well as appropriately value a reduced largest loss.
    \item We demonstrate the effectiveness and applicability of this marginal-pricing methodology through Economic Dispatch (ED) and Unit Commitment (UC) case studies. These studies demonstrate that the proposed pricing scheme for frequency services puts in place the right incentives: FR providers are incentivised to provide faster FR because they are rewarded according to their speed of delivery, and RES generators are incentivised to provide frequency services (such as synthetic inertia). 
\end{enumerate}

The rest of the paper is organised as follows: Section \ref{SectionDeductionConstraints} presents the post-fault frequency security constraints. These constraints are used in Section \ref{SectionPricing} to develop the marginal-pricing scheme for frequency services, while this pricing scheme is applied to several case studies in Section \ref{SectionCaseStudies}. Finally, Section \ref{SectionConclusions} gives the conclusion.

\section{Constraints for Frequency Security} \label{SectionDeductionConstraints}

Post-fault frequency dynamics are mathematically described by the swing equation \cite{KundurBook}:
\begin{equation} \label{SwingEqMultiFR}
\frac{2H}{f_0}\cdot\frac{\textrm{d} \Delta f(t)}{\textrm{d} t} = \textrm{FR}(t) - P_{\textrm{L}}
\end{equation}
Note that load damping has not been included eq. (\ref{SwingEqMultiFR}), which entails a conservative assumption used in this paper.

Graphically, the time-evolution of frequency deviation from nominal value following a large power outage takes this shape:
\begin{figure}[H]
\hspace*{0.05cm}
    \input{FrequencySimulation.tikz}
    \caption{Example of system frequency deviation after a power outage. The most extreme case still within frequency-security boundaries is shown, where the RoCoF attains the maximum value of $\textrm{RoCoF}_\textrm{max}$, and the frequency nadir drops to $\Delta f_\textrm{max}$.}
	\label{FigFrequencySimulation}
\end{figure}
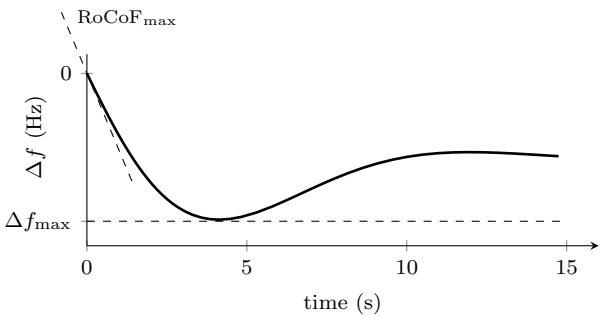

By solving the swing equation, the conditions for respecting the RoCoF, nadir and quasi-steady-state limits can be obtained. To illustrate the procedure for obtaining these conditions, consider the following case with three FR-services, where the aggregated system Frequency Response, $\textrm{FR}(t)$ in (\ref{SwingEqMultiFR}), is defined as:
\begin{subequations} \label{FRdefinition_3_NoDelay}
\begin{empheq}[left ={\textrm{FR}(t)\hspace{-2pt}=\hspace{-2pt}\empheqlbrace}]{alignat=3}
 & \left(\frac{R_1}{\textrm{T}_1}+\frac{R_2}{\textrm{T}_2}+\frac{R_3}{\textrm{T}_3}\right) \cdot t && \;\;\; \mbox{if $t\leq\textrm{T}_1$} \tag{\ref{FRdefinition_3_NoDelay}.1} \\
 & R_1 + \left(\frac{R_2}{\textrm{T}_2}+\frac{R_3}{\textrm{T}_3}\right) \cdot t && \;\;\; \mbox{if $\textrm{T}_1<t\leq\textrm{T}_2$} \tag{\ref{FRdefinition_3_NoDelay}.2}\\
 & R_1 + R_2 + \frac{R_3}{\textrm{T}_3} \cdot t && \;\;\; \mbox{if $\textrm{T}_2<t\leq\textrm{T}_3$} \tag{\ref{FRdefinition_3_NoDelay}.3}\\
 & R_1 + R_2 + R_3 && \;\;\; \mbox{if $t> \textrm{T}_3$} \tag{\ref{FRdefinition_3_NoDelay}.4} 
\end{empheq}
\end{subequations}
Function $\textrm{FR}(t)$ is graphically represented in Fig. \ref{Fig3services}:
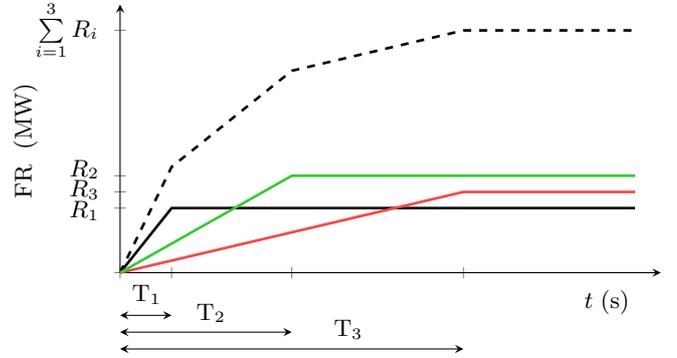
\begin{figure}[H]
\raggedright
\hspace*{-0.45cm}
    \input{FR_services.tikz}
    \caption{Time evolution of three distinct FR services that start ramping up at the very moment of the power outage, $t=0$. The total system $\textrm{FR}(t)$ as defined in eq. (\ref{FRdefinition_3_NoDelay}) is given by the dashed line.}
	\label{Fig3services}
\end{figure}

The different FR services in eq. (\ref{FRdefinition_3_NoDelay}) have been modelled using a linear ramp, while the actual dynamics of FR provided by generators are typically driven by droop controllers. This linear-ramp assumption for FR allows to deduce closed-form conditions for the frequency-security requirements presented in this section, and this assumption was demonstrated in \cite{OPFChavez} to approximate a real droop control while still maintaining system security.

The highest RoCoF occurs at $t=0$, when the power imbalance is highest because the power contribution from Frequency Response is still zero:
\begin{equation} 
\left|\mbox{RoCoF}(t=0)\right| = \frac{P_{\textrm{L}}\cdot f_0}{2 \cdot H} 
\end{equation}
Therefore, by limiting this highest RoCoF to be below the RoCoF requirement defined by the system operator, the following constraint for respecting RoCoF security is obtained:
\begin{equation} \label{RocofConstraint}
 \frac{P_{\textrm{L}}\cdot f_0}{\mbox{RoCoF}_{\textrm{max}}} \leq 2 \cdot H  
\end{equation}

In order for frequency to stabilize in a quasi-steady-state before Secondary and Tertiary Frequency Response take frequency back to the nominal 50Hz, the following q-s-s condition must be met:
\begin{equation} \label{qssMultiFR}
R_1 + R_2 + R_3 \geq P_{\textrm{L}}
\end{equation}
Which states that the total amount of FR must be at least equal to the value of the power outage $P_{\textrm{L}}$.

Finally, in order to respect all three frequency security requirements, the frequency nadir must be guaranteed to be above $\Delta f_{\textrm{max}}$. Before deducing the condition for nadir security, it is key to notice that the nadir will occur at a time-instant $t_\textrm{nadir}$ when the RoCoF becomes zero. The analytical expression for this $t_\textrm{nadir}$ can be obtained by solving the following equation, obtained by setting RoCoF to zero in the swing eq. (\ref{SwingEqMultiFR}):
\begin{equation} \label{PowerEquilibrium}
\textrm{FR}(t_\textrm{nadir}) = P_{\textrm{L}}
\end{equation}
Since three distinct FR services are considered for $\textrm{FR}(t)$ in (\ref{FRdefinition_3_NoDelay}), the time when the nadir occurs will depend on the total amount of FR delivered by each service, i.e. it will depend on the value of $R_1$, $R_2$ and $R_3$. As $R_1$, $R_2$ and $R_3$ will be decision variables in any optimisation problem which minimises the cost of frequency services, it is impossible to know exactly when the nadir will occur before solving the optimisation, and therefore every possibility for nadir must be considered beforehand. 

Let's start by assuming that the nadir happens in the first time-interval defined in (\ref{FRdefinition_3_NoDelay}), i.e. the time-interval when $t\leq\textrm{T}_1$. The power equilibrium defined by (\ref{PowerEquilibrium}) must then happen in that first time-interval, which means that the following condition holds:
\begin{equation} \label{conditionNadirInterval1}
R_1 + R_2 \frac{\textrm{T}_{1}}{\textrm{T}_{2}} + R_3 \frac{\textrm{T}_{1}}{\textrm{T}_{3}} > P_{\textrm{L}}
\end{equation}
The time within this first time-interval at which nadir is exactly reached is given by solving (\ref{PowerEquilibrium}) for that particular time-interval:
\begin{equation} \label{t_nadir_1}
t_\textrm{nadir} = \frac{P_{\textrm{L}}}{R_1/\textrm{T}_1+R_2/\textrm{T}_2+R_3/\textrm{T}_3} 
\end{equation}
One last step is necessary for obtaining the constraint which guarantees nadir security in this first time-interval $t \in [0,\textrm{T}_1]$, and it involves obtaining the expression for frequency deviation during that time-interval. This can be done by solving (\ref{SwingEqMultiFR}) for $t \in [0,\textrm{T}_1]$:
\begin{equation} \label{Sol_Delta_f_1}
\left|\Delta f(t)\right| = \frac{f_0}{2 H}\left(P_{\textrm{L}} \cdot t - \sum_{i=1}^{3}\frac{R_i}{2\textrm{T}_i} \cdot t^2  \right)
\end{equation}
Lastly, by substituting $t=t_\textrm{nadir}$ into (\ref{Sol_Delta_f_1}) using the expression in (\ref{t_nadir_1}), the value of frequency deviation at the nadir can be deduced:
\begin{equation} 
\left| \Delta f(t=t_\textrm{nadir}) \right| =
\frac{f_0}{4 H} \cdot \frac{(P_\textrm{L})^2}{R_1/\textrm{T}_1+R_2/\textrm{T}_2+R_3/\textrm{T}_3}
\end{equation}
In conclusion, by limiting this maximum frequency deviation at the nadir to be below the requirement established by the system operator, i.e. enforcing $\left| \Delta f(t=t_\textrm{nadir})  \right| \leq \Delta f_{\textrm{max}}$, the constraint for respecting nadir security is obtained. This constraint takes the form:
\begin{equation} \label{nadir_req_1}
\frac{H}{f_0} \cdot  \left(\frac{R_1}{\textrm{T}_{1}} + \frac{R_2}{\textrm{T}_{2}} + \frac{R_3}{\textrm{T}_{3}} \right) \geq 
\frac{(P_{\textrm{L}})^2}{4 \Delta f_{\textrm{max}}}
\end{equation}

Remember that in order to obtain (\ref{nadir_req_1}) we have assumed that the nadir happens during the first time interval for $\textrm{FR}(t)$, which means that condition (\ref{conditionNadirInterval1}) holds. However, the nadir could also occur during any other time-interval defined in (\ref{FRdefinition_3_NoDelay}): time-interval 2 where $t \in (\textrm{T}_1,\textrm{T}_2]$ or time-interval 3 where $t \in (\textrm{T}_2,\textrm{T}_3]$. The nadir constraints that would apply if the nadir happened in any of these time-intervals can be deduced by following the same mathematical procedure described above for time-interval 1. Therefore, all following constraints, along with the conditions for each of them to be enforced, must be implemented in an optimisation problem to guarantee nadir security:

if $\quad R_1 + R_2 \frac{\textrm{T}_{1}}{\textrm{T}_{2}} + R_3 \frac{\textrm{T}_{1}}{\textrm{T}_{3}} > P_{\textrm{L}} \quad$ then enforce:
\begin{equation} \label{nadirThreeServices1}
\hspace*{14mm}
\frac{H}{f_0} \cdot  \left(\frac{R_1}{\textrm{T}_{1}} + \frac{R_2}{\textrm{T}_{2}} + \frac{R_3}{\textrm{T}_{3}} \right) \geq 
\frac{(P_{\textrm{L}})^2}{4 \Delta f_{\textrm{max}}}
\end{equation}

else if $\quad R_1 + R_2 + R_3 \frac{\textrm{T}_{2}}{\textrm{T}_{3}} > P_{\textrm{L}} \quad$ then enforce:
\begin{multline} \label{nadirThreeServices2}
\vspace*{-5mm}
\left(\frac{H}{f_0} - \frac{R_1\textrm{T}_{1}}{4 \Delta f_{\textrm{max}}}\right) \cdot  \left(\frac{R_2}{\textrm{T}_{2}} + \frac{R_3}{\textrm{T}_{3}} \right) \geq 
\frac{(P_{\textrm{L}}-R_1)^2}{4 \Delta f_{\textrm{max}}}
\end{multline}

else enforce:
\begin{multline} \label{nadirThreeServices3}
\left(\frac{H}{f_0} - \frac{R_1\textrm{T}_{1}+R_2\textrm{T}_{2}}{4 \Delta f_{\textrm{max}}}\right) \cdot  \frac{R_3}{\textrm{T}_{3}} \geq 
\frac{(P_{\textrm{L}}-R_1-R_2)^2}{4 \Delta f_{\textrm{max}}}
\end{multline}

The first condition states that the power equilibrium happens before $\textrm{T}_1$, when none of the three FR services have finished ramping up by the nadir; the second condition states that the power equilibrium happens after $\textrm{T}_1$ but before $\textrm{T}_2$, then only $\textrm{FR}_{1}$ has finished ramping up by the nadir; finally, the third condition states that the power equilibrium happens after $\textrm{T}_2$, then both $\textrm{FR}_{1}$ and $\textrm{FR}_{2}$ have finished ramping up (note that the nadir must always take place before the slowest FR service finishes ramping up, as otherwise it would be impossible to reach a power equilibrium and frequency would drop indefinitely). These conditional constraints (\ref{nadirThreeServices1}), (\ref{nadirThreeServices2}) and (\ref{nadirThreeServices3}), which are only enforced if their associated condition is met, can be implemented in an optimisation problem with the use of auxiliary binary variables, refer to \cite{ConejoOptimizationBook} and \cite{YALMIPconditional} for further details.

In conclusion: the RoCoF and q-s-s constraints, (\ref{RocofConstraint}) and (\ref{qssMultiFR}), along with the conditional nadir constraints (\ref{nadirThreeServices1}), (\ref{nadirThreeServices2}) and (\ref{nadirThreeServices3}), guarantee that system frequency will stay within the margins defined by the system operator in the event of a power outage. The feasible space defined by these five constraints is an MISOCP: (\ref{RocofConstraint}) and (\ref{qssMultiFR}) are linear constraint, while (\ref{nadirThreeServices1}), (\ref{nadirThreeServices2}) and (\ref{nadirThreeServices3}) are rotated SOCs. The binary variables needed to implement the conditions in the nadir constraints make the problem an MISOCP.

Note that a particular case of three FR services with distinct delivery times has been considered in order to deduce the frequency-security conditions in this section. Nevertheless, any number of frequency services could be considered, and following the same mathematical procedure included in this section, MISOCP constraints would be obtained.

\subsection{Considering time-delays in FR provision} \label{SectionDelays}

Here we demonstrate that the frequency-security constraints deduced in the previous section can also be obtained in MISOCP form for a more general case: the case in which FR services might be activated some time after the power outage, i.e. they have an activation time-delay before starting ramping up. This activation delay would mean that these FR services do not react to every tiny frequency deviation from nominal state, therefore reducing wear and tear for thermal units and greatly decreasing the number of times that battery storage devices providing FR switch from charging to discharging mode (and vice versa).

Consider the case of two distinct frequency services with also distinct activation delays:
\begin{subequations} \label{FRdefinition_2_delays}
\begin{empheq}[left ={\textrm{FR}(t)\hspace{-2pt}=\hspace{-2pt}\empheqlbrace}]{alignat=3}
 & \qquad 0 && \;\; \mbox{if $t\leq\textrm{T}_{\textrm{del},1}$} \tag{\ref{FRdefinition_2_delays}.1} \\[-2pt]
 & \frac{R_1}{\textrm{T}_1} (t-\textrm{T}_{\textrm{del},1}) && \;\; \mbox{if $\textrm{T}_{\textrm{del},1}<t\leq\textrm{T}_1$} \tag{\ref{FRdefinition_2_delays}.2}\\[-2pt]
 & \qquad R_1 && \;\; \mbox{if $\textrm{T}_1<t\leq \textrm{T}_{\textrm{del},2}$} \tag{\ref{FRdefinition_2_delays}.3}\\[1pt]
 & R_1 + \frac{R_2}{\textrm{T}_2} (t-\textrm{T}_{\textrm{del},2}) && \;\; \mbox{if $\textrm{T}_{\textrm{del},2}<t\leq \textrm{T}_2$} \tag{\ref{FRdefinition_2_delays}.4} \\[1pt]
 & \quad R_1 + R_2 && \;\; \mbox{if $t> \textrm{T}_2$} \tag{\ref{FRdefinition_2_delays}.5}
\end{empheq}
\end{subequations}
Graphically, these FR services take this shape:
\begin{figure}[H]
\raggedright
\hspace*{-0.45cm}
    \input{FR_services2.tikz}
    \caption{Time evolution of two distinct FR services, each one with a different activation time-delay. The total system $\textrm{FR}(t)$ as defined in eq. (\ref{FRdefinition_2_delays}) is given by the dashed line.}
	\label{Fig2services_delay}
\end{figure}
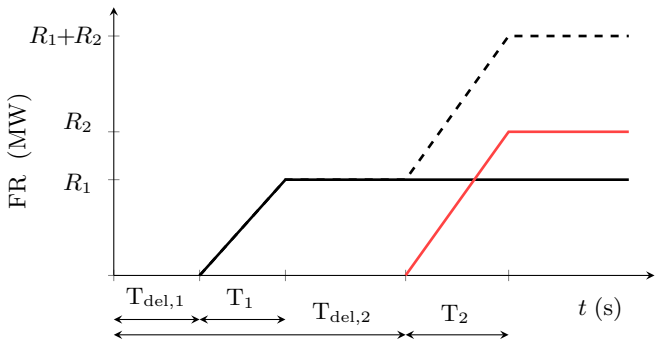

The RoCoF and q-s-s constraints for this case are still (\ref{RocofConstraint}) and (\ref{qssMultiFR}). The nadir constraints can be deduced following the same procedure as in Section \ref{SectionDeductionConstraints}: obtain the analytical expression for $t_\textrm{nadir}$ for each of the time intervals defined in (\ref{FRdefinition_2_delays}), solve the swing equation (\ref{SwingEqMultiFR}) for each of these time-intervals, and finally substitute the expression for $t_\textrm{nadir}$ in the just obtained solution for swing equation. Following these steps, the conditional nadir constraints for this case are:

if $\quad R_1 > P_{\textrm{L}} \quad$ then enforce:
\begin{multline} \label{nadirTwoServices_delays1}
\left(\frac{H}{f_0} + \frac{R_1 \cdot \textrm{T}_{\textrm{del},1}^2/\textrm{T}_{1}}{4 \Delta f_{\textrm{max}}}\right) \cdot  \frac{R_1}{\textrm{T}_{1}}  \geq 
\frac{(P_{\textrm{L}} + R_1 \cdot \textrm{T}_{\textrm{del},1}/\textrm{T}_{1})^2}{4 \Delta f_{\textrm{max}}}
\end{multline}

\vspace{-2mm}
else enforce:
\begin{multline} \label{nadirTwoServices_delays2}
\vspace*{-5mm}
\left(\frac{H}{f_0} - \frac{R_1 (\textrm{T}_{1}\hspace{-1mm}+\hspace{-1mm}2\textrm{T}_{\textrm{del},1})}{4 \Delta f_{\textrm{max}}}  + \frac{R_2 \cdot \textrm{T}_{\textrm{del},2}^2/\textrm{T}_{2}}{4 \Delta f_{\textrm{max}}} \right) \cdot  \frac{R_2}{\textrm{T}_{2}}  \\ 
\geq 
\frac{(P_{\textrm{L}} - R_1 + R_2 \cdot \textrm{T}_{\textrm{del},2}/\textrm{T}_{2})^2}{4 \Delta f_{\textrm{max}}}
\end{multline}
Note that nadir will never happen at time-interval 1 in (\ref{FRdefinition_2_delays}) because there is no FR delivered during that interval, so it is impossible to reach a power equilibrium. Nadir will also not happen in time-interval 3, because function $\textrm{FR}(t)$ is constant during that time-interval. Therefore only time-intervals 2 and 4 in (\ref{FRdefinition_2_delays}) have been considered for deducing (\ref{nadirTwoServices_delays1}) and (\ref{nadirTwoServices_delays2}).

Both constraints (\ref{nadirTwoServices_delays1}) and (\ref{nadirTwoServices_delays2}) are SOCs, and therefore the frequency-secured optimisation problem is still and MISOCP, due to the binary variables introduced by the conditional statements that enforce either (\ref{nadirTwoServices_delays1}) or (\ref{nadirTwoServices_delays2}).

Finally, it is worth mentioning that although a particular case of just two FR services with activation delays has been considered for deducing (\ref{nadirTwoServices_delays1}) and (\ref{nadirTwoServices_delays2}), any combination of FR services with distinct delays would still yield nadir constraints in SOC form.

\section{Marginal-Pricing Scheme for Frequency Services} \label{SectionPricing}

In this section we propose a settlement mechanism which assigns shadow prices to the different frequency services so that an efficient market equilibrium can be obtained. This settlement is based on maximization of social welfare using strong duality. 

The optimisation problem defined by the frequency-security constraints contains integer variables, and therefore strong duality does not hold. The sources of integrality are two: 1) $H$, since inertia is a discrete magnitude dependent on the on/off state of thermal generators; and 2) Binary variables introduced by the conditional statements for the nadir constraints, as in (\ref{nadirTwoServices_delays1}) and (\ref{nadirTwoServices_delays2}). Here we propose a methodology to overcome these issues and allow to define prices for the different frequency services.

Regarding inertia, the approach proposed in \cite{ElaI} is adopted, consisting on relaxing the integer decision variable to a continuous one. Although this approach would not lead to a feasible dispatch, a parallel optimisation maintaining the integrality of $H$ would be run, giving the optimal dispatch. The relaxed problem would then give the market settlement allowing to price inertia as a service.

For the auxiliary binary variables in the nadir conditional statements, the following approach is proposed: first, the original MISOCP is run, which chooses the optimal time for nadir to occur by enforcing only one of the nadir constraints associated with a particular time-interval; then, a new optimisation is run, including only the nadir constraint that was chosen in the previous step and removing any binary variables, therefore yielding a convex SOCP formulation. 
Using this approach, strong duality holds for the resulting convex problem.

The price of each frequency service can be obtained from the Karush Kuhn Tucker (KKT) conditions of the frequency-constrained optimisation problem. The price-term for each service corresponding to the RoCoF (\ref{RocofConstraint}) and q-s-s (\ref{qssMultiFR}) constraints is immediately obtained from the dual variable of each constraint, since they are linear. These constraints appear in the Lagrangian with their corresponding Lagrange multipliers: 

\begin{equation} 
\lambda_\textrm{R} \left(\frac{P_{\textrm{L}}\cdot f_0}{\mbox{RoCoF}_{\textrm{max}}} - 2 \cdot H\right)
\end{equation}

\begin{equation} 
\lambda_\textrm{qss}[P_{\textrm{L}} - (R_1 + R_2)]
\end{equation}

For the deduction of the price-terms associated with the nadir constraint in this section, the two FR services with activation delays considered in section \ref{SectionDelays} are also considered. This is because that case allows to present the most general pricing scheme needed, corresponding to nadir happening in time-interval 4 of (\ref{FRdefinition_2_delays}): the nadir constraint enforced is then (\ref{nadirTwoServices_delays2}), for which service $\textrm{FR}_1$ has finished ramping up while service $\textrm{FR}_2$ is still ramping up. Therefore, any service that has finished ramping up by the nadir can be priced in the same way as $\textrm{FR}_1$ in this example, while any service still ramping up can be priced as $\textrm{FR}_2$.

Before deducing the pricing terms associated with constraint (\ref{nadirTwoServices_delays2}), it must first be re-formulated to be expressed in standard SOC form \cite{ConvexOptPowerBook}:
\begin{multline}
\norm{
\begin{bmatrix} \frac{1}{f_0}  &\frac{-\textrm{T}_1-\textrm{T}_{\textrm{del},1}}{4\Delta f_\textrm{max}}  &\frac{\textrm{T}_{\textrm{del},2}^2/\textrm{T}_{2}}{4\Delta f_\textrm{max}}-\frac{1}{\textrm{T}_{2}} &0 \\[5pt]
0 &\frac{-1}{\sqrt{\Delta f_\textrm{max}}} &\frac{\textrm{T}_{\textrm{del},2}/\textrm{T}_2}{\sqrt{\Delta f_\textrm{max}}} &\frac{1}{\sqrt{\Delta f_\textrm{max}}}
\end{bmatrix}
\begin{bmatrix} H  \\ R_1 \\ R_2 \\ P_\textrm{L} \end{bmatrix} 
} 
\\
\leq  \begin{bmatrix} \frac{1}{f_0}  &\frac{-\textrm{T}_1-\textrm{T}_{\textrm{del},1}}{4\Delta f_\textrm{max}}  &\frac{\textrm{T}_{\textrm{del},2}^2/\textrm{T}_{2}}{4\Delta f_\textrm{max}}+\frac{1}{\textrm{T}_{2}} &0 
\end{bmatrix}
\begin{bmatrix} H  \\ R_1 \\ R_2 \\ P_\textrm{L} \end{bmatrix} 
\end{multline}

This constraint appears in the Lagrangian with its respective vector multipliers:
\begin{multline}
\lambda_1 \Biggl(\frac{H}{f_0} -R_1 \frac{\textrm{T}_1+\textrm{T}_{\textrm{del},1}}{4\Delta f_\textrm{max}}  +R_2\biggl[\frac{\textrm{T}_{\textrm{del},2}^2/\textrm{T}_{2}}{4\Delta f_\textrm{max}}-\frac{1}{\textrm{T}_{2}}\biggr]\Biggr) 
\\
+ \lambda_2 \Biggl(\frac{P_\textrm{L}-R_1}{\sqrt{\Delta f_\textrm{max}}}  +R_2\frac{\textrm{T}_{\textrm{del},2}/\textrm{T}_{2}}{\sqrt{\Delta f_\textrm{max}}}\Biggr)
\\
\hspace*{-2mm}- \mu \Biggl(\frac{H}{f_0} -R_1 \frac{\textrm{T}_1+\textrm{T}_{\textrm{del},1}}{4\Delta f_\textrm{max}}  +R_2\biggl[\frac{\textrm{T}_{\textrm{del},2}^2/\textrm{T}_{2}}{4\Delta f_\textrm{max}}+\frac{1}{\textrm{T}_{2}}\biggr]\Biggr)
\end{multline}

The following constraint must be included in the dual problem to comply with the KKT condition for dual feasibility, due to using vector dual variables for the SOC:
\begin{equation}
\norm{\begin{matrix}\lambda_1 \\ \lambda_2 \end{matrix}} \leq \mu
\end{equation}

The KKT stationarity condition is given by the gradient of the Lagrangian. Then, by differentiating the Lagrangian with respect to each frequency service, the corresponding price is obtained:
\begin{itemize}
    \item The price of $H$ would be: 
    \begin{equation}
    \underbrace{\frac{\mu-\lambda_1}{f_0}}_\text{nadir} 
    \underbrace{\vphantom{\frac{\mu-\lambda_1}{f_0}} + 2 \lambda_\textrm{R}}_\text{RoCoF}
    \end{equation}
    
    \item The price of $P_\textrm{L}$ is:
    \begin{equation}
    \underbrace{-\frac{\lambda_2}{\sqrt{\Delta f_\textrm{max}}}}_\text{nadir}  
    \underbrace{\vphantom{\frac{\lambda_2}{\sqrt{\Delta f_\textrm{max}}}} - \frac{\lambda_\textrm{R}\cdot f_0}{\mbox{RoCoF}_{\textrm{max}}} }_\text{RoCoF}
    \underbrace{\vphantom{\frac{\lambda_2}{\sqrt{\Delta f_\textrm{max}}}} - \lambda_\textrm{qss}}_\text{q-s-s}
    \end{equation}
    Note that the price for $P_\textrm{L}$ would be negative since a larger loss would increase the cost for ancillary-service provision.
    
    \item The price of an FR service which has finished ramping up by nadir (as $\textrm{FR}_1$ in this example) is given by:
    \begin{equation}
    \underbrace{\frac{\lambda_2}{\sqrt{\Delta f_\textrm{max}}}
    -(\mu-\lambda_1) \frac{\textrm{T}_1+\textrm{T}_{\textrm{del},1}}{4\Delta f_\textrm{max}}}_\text{nadir} 
    \underbrace{\vphantom{\frac{\lambda_2}{\sqrt{\Delta f_\textrm{max}}}} + \lambda_\textrm{qss}}_\text{q-s-s}
    \end{equation}
    Therefore the value of the service decreases with increasing delivery time and increasing delay. In other words, the service becomes less valuable as its dynamics are slower.
    
    \item The price of FR services still ramping up by the nadir (as $\textrm{FR}_2$ in this example) is given by:
    \begin{multline}
    \underbrace{(\mu+\lambda_1)\frac{1}{\textrm{T}_{2}}
    +
    (\mu-\lambda_1)\frac{\textrm{T}_{\textrm{del},2}^2/\textrm{T}_{2}}{4\Delta f_\textrm{max}}
    - \lambda_2 \frac{\textrm{T}_{\textrm{del},2}/\textrm{T}_2}{\sqrt{\Delta f_\textrm{max}}}}_\text{nadir} \\
    \underbrace{\vphantom{\frac{\lambda_2}{\sqrt{\Delta f_\textrm{max}}}} + \lambda_\textrm{qss}}_\text{q-s-s}
    \end{multline}
\end{itemize}

The pricing methodology proposed in this section applies to any convex optimisation problem which includes the frequency-security constraints deduced in Section \ref{SectionDeductionConstraints}, including the case in which the energy and FR markets are cleared simultaneously.

\subsection{Discussion on uplift}

The inherent non-convexities in power systems, where the commitment decision of thermal units is forcibly a binary variable, makes it impossible to apply duality theory even to price energy in electricity markets. This issue has been studied in the literature and several approaches have been proposed to overcome this challenge and obtain marginal prices for energy \cite{ConvexHullPricing}:
\begin{itemize}
    \item Dispatchable model, where the commitment decision of thermal units is relaxed from binary to continuous, as used in the present paper.

    \item Restricted model, which consist on first running the original mixed-integer problem, then running a second optimisation relaxing the commitment decision but using equality constraints to force the commitment variables to take the same value as in the solution of the original, mixed-integer problem. Some studies have observed that this model can lead to high volatility of energy prices \cite{ConvexOptPowerBook}.

    \item Convex hull approach, which consists on substituting the original mixed-integer problem by its convex hull, i.e. the tightest possible convex relaxation of the original problem. The authors in \cite{ConvexHullPricing} argue that this approach provides the smallest make-whole payments, while computing the convex hull of a problem is a very computationally-intensive task.
\end{itemize}
This is still an active area of research, as no method has been proven to be superior in all instances for the pricing of ``uplift" (some very recent publications on the matter are available \cite{Hobbs2019Pricing}). The term ``uplift" refers to the payments that must be made to certain generators in order to compensate for the negative revenue that they might obtain at times, given that the energy prices only support the optimal solution in the simplified, convex model where they are derived from, but not in the original non-convex model that more accurately represents the real power system.

Since the pricing of uplift is inherently related to the non-convexities in a power system and not just related to inertia, it is out of the scope of the present paper to study uplift in detail. Therefore, we have opted to use the dispatchable model for allowing to price inertia, as was proposed in \cite{ElaI}. However, our pricing formulation could also be applied if any of the other methods for pricing uplift is considered. Future work should study the implications of other uplift-pricing models for the co-optimisation of energy and frequency services.

\section{Case Studies} \label{SectionCaseStudies}

In this section we carry out several frequency-secured ED and UC optimisations, to demonstrate the applicability and implications of our proposed marginal-pricing methodology for frequency services. These optimisation problems, formulated as MISOCPs due to the conditional nadir constraints, were implemented in YALMIP \cite{YALMIP}, an optimisation modelling layer for MATLAB, while the low-level numerical computations were solved by calling commercial solver Gurobi. The code used for solving these examples is available in \href{https://github.com/badber/MultiFR_optimisation}{\underline{this link}}. 

In these studies, the different generators submit a bid for providing energy to the market operator, and the market clearing consists on minimising total system costs. There is no explicit bid for inertia or FR provision. The frequency-security requirements are set to $\textrm{RoCoF}_{\textrm{max}}=1\textrm{Hz/s}$ and $\Delta f_{\textrm{max}}=0.8\textrm{Hz}$.

\subsection{Value of single-speed FR} \label{CaseStudySingleFR}

This section demonstrates that if there are no distinct speeds for FR, the pricing scheme for the single FR service is equivalent to that of reserve: generators simply recover their opportunity cost. This study is based on Example 6.7 in \cite{GoranBook}, which demonstrates that by using the dual variable from the reserve constraint in a simultaneous clearing of the energy and reserve markets, all the generators would recover the opportunity cost from limiting their power output due to providing reserve.

To illustrate this point, we run a single-time-period ED, simultaneously clearing the energy and FR markets. The optimisation problem is formulated as:

\begin{alignat}{3}
& {\text{min}} \quad
&& \sum_{g \in \mathcal{G}}\textrm{c}^\textrm{m}_g \cdot P_g \label{ObjectiveFunction}\\
& \text{s.t.} \quad
&& \sum_{g \in \mathcal{G}}P_g = \textrm{P}_\textrm{D} \label{LoadBalance}\\
& && 0 \leq P_g \leq \textrm{P}^\textrm{max}_g && \quad \forall g \in \mathcal{G} \label{GenLimits}\\
& && 0 \leq R_g \leq \textrm{R}^\textrm{max}_g && \quad \forall g \in \mathcal{G} \label{FR_Limits}\\
& && 0 \leq R_g \leq \textrm{P}^\textrm{max}_g-P_g && \quad \forall g \in \mathcal{G} \label{FR_Limits2}\\
& && \textrm{RoCoF} \; \textrm{constraint} \nonumber\\
& && \textrm{Nadir} \; \textrm{constraints} \nonumber\\
& && \textrm{q-s-s} \; \textrm{constraint} \nonumber
\end{alignat}
Constraint (\ref{LoadBalance}) enforces the power balance, (\ref{GenLimits}) enforce generation limits, (\ref{FR_Limits}) enforce FR-provision limits, and (\ref{FR_Limits2}) limits FR provision to be below the headroom of each generator.

The characteristics of generators participating in the energy and FR markets are given in Table \ref{TableThermal_ED}.
\begin{table}[H]
\renewcommand{\arraystretch}{1.2}
\small
\centering
\begin{tabular}{p{3.3cm}| l| p{1.2cm}| p{1.2cm}}
    \multicolumn{1}{c|}{} & Nuclear & Gen. type 1 & Gen. type 2\\
\hline
Power range (MW) & 100-100 & 0-400 & 0-300\\ 
Number of Units & 1 & 5 & 5\\ 
FR capacity (MW) & 0 & 225 & 175\\ 
Energy bid (\pounds/MWh) & 15 & 17 & 18\\ 
Inertia constant (s) & 6 & 6 & 6\\ 
FR delivery time (s) & N/A & 10 & 10\\ \hline
\end{tabular}
\caption{\label{TableThermal_ED}Characteristics of thermal generators}
\end{table}
All generators willing to provide FR are assumed to provide Primary Frequency Response with a delivery time of 10s, which was the only FR service defined by National Grid in Great Britain until 2016. The largest possible outage is driven by the largest single plant, the nuclear unit, therefore setting $P_\textrm{L}=\textrm{P}_{\textrm{L}}^{\textrm{max}}=100\textrm{MW}$. Note that there are 5 units of Generators type 1, therefore each of them has a rating of 80MW. Similarly, 5 units of Generators type 2 are available, with a 60MW-rating each. 

Some further clarifications on Table \ref{TableThermal_ED}: for the sake of simplicity, we choose zero Minimum Stable Generation (MSG) for all generators, to better illustrate the pricing of FR services (zero MSG was also used in Example 6.7 in \cite{GoranBook}). Also for clarity, we consider fixed system inertia; however, the co-optimisation of inertia and FR is possible with our model, and will be demonstrated in a UC example in Section \ref{SectionInertia}.

Since both types of generators 1 and 2 in Table \ref{TableThermal_ED} have the same FR delivery-time of 10s, the nadir constraint is: 
\begin{equation} \label{nadir_first_ED}
R_1 \geq 
\frac{(P_{\textrm{L}})^2\cdot \textrm{T}_1 \cdot f_0}{4 \Delta f_{\textrm{max}}\cdot H}
= 372\textrm{MW}
\end{equation}
Which is equivalent to a fixed reserve requirement, such as $\textrm{Reserve}\geq 372\textrm{MW}$.

The results for a low-demand case of $\textrm{P}_\textrm{D}=250\textrm{MW}$ are presented in Table \ref{TableResults_ED_LowD}:
\begin{table}[H]
\renewcommand{\arraystretch}{1.2}
\small
\centering
\begin{tabular}{p{3.6cm}| l| p{1.1cm}| p{1.1cm}}
    \multicolumn{1}{c|}{} & Nuclear & Gen. type 1 & Gen. type 2\\
\hline
Power produced (MW) & 100 & 150 & 0\\ 
FR provided (MW) & 0 & 197 & 175\\ 
Revenue from energy (\pounds) & 1700 & 2550 & 0\\
Revenue from FR (\pounds) & 0 & 0 & 0\\ 
Profit (\pounds) & 200 & 0 & 0\\ \hline
Energy price (\pounds/MWh) & \multicolumn{3}{c}{17}\\ \hline
FR price (\pounds/MW) & \multicolumn{3}{c}{0}\\ \hline
\end{tabular}
\caption{\label{TableResults_ED_LowD}Solution of the frequency-secured ED for the system in Table \ref{TableThermal_ED}, for a low-demand case of $\textrm{P}_\textrm{D}=250\textrm{MW}$.}
\end{table}

The energy price has been calculated using the dual variable of the power-balance constraint (\ref{LoadBalance}), and the FR price comes from the dual variable of (\ref{nadir_first_ED}). Profit has been calculated assuming a competitive market in which energy bids are equal to marginal costs.

In this low-demand case, since FR is provided for free and all the energy can be provided by the cheapest Generators type 1, the price for FR is of \pounds0/MW. However, if the demand increases to reach $\textrm{P}_\textrm{D}=400\textrm{MW}$, the clearing prices change:

\begin{table}[H]
\renewcommand{\arraystretch}{1.2}
\small
\centering
\begin{tabular}{p{3.6cm}| l| p{1.1cm}| p{1.1cm}}
    \multicolumn{1}{c|}{} & Nuclear & Gen. type 1 & Gen. type 2\\
\hline
Power produced (MW) & 100 & 203 & 97\\ 
FR provided (MW) & 0 & 197 & 175\\ 
Revenue from energy (\pounds) & 1800 & 3654 & 1746\\
Revenue from FR (\pounds) & 0 & 197 & 175\\ 
Profit (\pounds) & 300 & 400 & 175\\ \hline
Energy price (\pounds/MWh) & \multicolumn{3}{c}{18}\\ \hline
FR price (\pounds/MW) & \multicolumn{3}{c}{1}\\ \hline
\end{tabular}
\caption{\label{TableResults_ED_HighD}Solution of the frequency-secured ED for the system in Table \ref{TableThermal_ED}, for a high-demand case of $\textrm{P}_\textrm{D}=400\textrm{MW}$.}
\end{table}

Now that the demand has increased, Generators type 2 start providing energy, so that Generators type 1 can still provide 197MW of FR (because the FR provision from Generators type 2 is capped at 175MW). The price of FR becomes \pounds1/MW while the price of energy is now of \pounds18/MW, so Generators type 1 are indifferent to providing FR or selling more energy: if they are part-loaded to provide FR, they recover the opportunity cost of not providing as much energy as they could. 

This study has demonstrated that if a single FR service is defined, as was common practice in Great Britain until recently, the same principles as for reserve pricing could be applied, since there was no significant difference between the single-speed FR and reserve. However, the following section will demonstrate the need for a marginal-pricing scheme tailored for FR if multi-speed FR services are to be recognised by the operator.

\subsection{Value and incentives for multi-speed FR} \label{CaseStudyMultiFR}

The case study presented here will demonstrate the difference between marginal-pricing for reserve and marginal-pricing for FR when more than one FR service is available.

Now let's consider the same system as in Section \ref{CaseStudySingleFR}, but assuming that Generators Type 1 can actually provide FR in 7s, and therefore they increase their energy bid to \pounds19/MWh due to providing this faster FR service (they then become the most expensive generators):
\begin{table}[H]
\renewcommand{\arraystretch}{1.2}
\small
\centering
\begin{tabular}{p{3.3cm}| l| p{1.2cm}| p{1.2cm}}
    \multicolumn{1}{c|}{} & Nuclear & Gen. type 1 & Gen. type 2\\
\hline
Power range (MW) & 100-100 & 0-400 & 0-300\\ 
Number of Units & 1 & 5 & 5\\ 
FR capacity (MW) & 0 & 225 & 175\\ 
Energy bid (\pounds/MWh) & 15 & \textbf{19} & 18\\ 
Inertia constant (s) & 6 & 6 & 6\\ 
FR delivery time (s) & N/A & \textbf{7} & 10\\ \hline
\end{tabular}
\caption{\label{TableThermal_ED_FR2}Characteristics of thermal generators. The only different parameters with respect to Table \ref{TableThermal_ED} are given in bold font.}
\end{table}

The results for a demand of $\textrm{P}_\textrm{D}=400\textrm{MW}$ are presented in Table \ref{TableResults_ED_HighD_FR2}:
\begin{table}[H]
\renewcommand{\arraystretch}{1.2}
\small
\centering
\begin{tabular}{p{3.6cm}| l| p{1.1cm}| p{1.1cm}}
    \multicolumn{1}{c|}{} & Nuclear & Gen. type 1 & Gen. type 2\\
\hline
Power produced (MW) & 100 & 50.6 & 249.4\\ 
FR provided (MW) & 0 & 225 & 50.6\\ 
Revenue from energy (\pounds) & 1900 & 961.4 & 4738.6\\
Revenue from FR (\pounds) & 0 & 315 & 50.6\\ 
Profit (\pounds) & 400 & 315 & 300\\ \hline
Energy price (\pounds/MWh) & \multicolumn{3}{c}{19}\\ \hline
$\textrm{FR}_1$ price (\pounds/MW) & \multicolumn{3}{c}{1.4}\\ \hline
$\textrm{FR}_2$ price (\pounds/MW) & \multicolumn{3}{c}{1}\\ \hline
\end{tabular}
\caption{\label{TableResults_ED_HighD_FR2}Solution of the frequency-secured ED for the system in Table \ref{TableThermal_ED_FR2}, for a demand of $\textrm{P}_\textrm{D}=400\textrm{MW}$.}
\end{table}

After solving the ED optimisation, one can verify that the enforced nadir constraint in this case is:
\begin{equation} 
\frac{R_1}{\textrm{T}_1} + \frac{R_2}{\textrm{T}_2} \geq 
\frac{(P_{\textrm{L}})^2\cdot f_0}{4 \Delta f_{\textrm{max}}\cdot H}
\end{equation}
Simplifying the above constraint:
\begin{equation} 
1.4 \cdot R_1 + R_2 \geq 372\textrm{MW}
\end{equation}

These results show that Generators type 2 are compensated for their part-loading by the payment for $\textrm{FR}_2$, therefore they recover their opportunity cost. On the other hand, Generators type 1 increase their revenue from providing the faster service $\textrm{FR}_1$ than if they were fully loaded to sell more energy. This demonstrates that fast provision of FR is incentivised under the proposed marginal-pricing scheme: FR providers can make a profit from providing faster FR, not simply recover the opportunity cost. Otherwise, there would be no incentive for FR providers to achieve faster delivery of FR. This is the key difference between reserve and FR: the speed of delivery of FR matters (while the dynamics of reserve do not matter), and that speed is reflected in the marginal pricing for FR that we propose. Therefore our methodology puts the right incentives in place to reduce whole system costs, as generators are incentivised to provide faster FR and that in turn benefits the whole system.

\subsection{Incentive for reducing FR activation delays} \label{CaseStudyDelays}

Here we study the impact of an activation delay in the value of an FR service, for which we use the same system as in Table \ref{TableThermal_ED_FR2} but considering that service $\textrm{FR}_1$ now has a 0.4 delay, i.e. $\textrm{T}_{\textrm{del,}1}=0.4\textrm{s}$.

The results for a demand of $\textrm{P}_\textrm{D}=400\textrm{MW}$ are presented in Table \ref{TableResults_ED_delay}:
\begin{table}[H]
\renewcommand{\arraystretch}{1.2}
\small
\centering
\begin{tabular}{p{3.6cm}| l| p{1.1cm}| p{1.1cm}}
    \multicolumn{1}{c|}{} & Nuclear & Gen. type 1 & Gen. type 2\\
\hline
Power produced (MW) & 100 & 143.5 & 156.5\\ 
FR provided (MW) & 0 & 225 & 143.5\\ 
Revenue from energy (\pounds) & 1900 & 2726.5 & 2973.5\\
Revenue from FR (\pounds) & 0 & 222.7 & 143.5\\ 
Profit (\pounds) & 400 & 222.7 & 300\\ \hline
Energy price (\pounds/MWh) & \multicolumn{3}{c}{19}\\ \hline
$\textrm{FR}_1$ price (\pounds/MW) & \multicolumn{3}{c}{0.99}\\ \hline
$\textrm{FR}_2$ price (\pounds/MW) & \multicolumn{3}{c}{1}\\ \hline
\end{tabular}
\caption{\label{TableResults_ED_delay}Solution of the frequency-secured ED for the system in Table \ref{TableThermal_ED_FR2}, for a demand of $\textrm{P}_\textrm{D}=400\textrm{MW}$. Service $\textrm{FR}_1$ has a 0.4s activation delay in this case.}
\end{table}

The enforced nadir constraint in this case is:
\begin{multline} 
\left(\frac{H}{f_0} + \frac{R_1 \cdot \textrm{T}_{\textrm{del},1}^2/\textrm{T}_{1}}{4 \Delta f_{\textrm{max}}} \right)\cdot \left(\frac{R_1}{\textrm{T}_1} + \frac{R_2}{\textrm{T}_2}\right) \geq \\
\frac{(P_{\textrm{L}}+R_1\cdot\textrm{T}_{\textrm{del,}1}/\textrm{T}_1)^2}{4 \Delta f_{\textrm{max}}}
\end{multline}

The price of service $\textrm{FR}_1$ has dropped to \pounds0.99/MW, making $\textrm{FR}_1$ less valuable than $\textrm{FR}_2$ although the delivery time of $\textrm{FR}_2$ is still longer than $\textrm{T}_1$, but $\textrm{FR}_2$ still has no delay. Due to the 0.4s activation delay, the profit of Generators type 1 has decreased by 29\% compared to the results in Table \ref{TableResults_ED_HighD_FR2}, which demonstrates how FR providers are incentivised to reduce their activation delay with our marginal-pricing methodology.

\subsection{Value of a reduced largest-loss} \label{SectionPLoss}

Now we analyse the value of a “reduced largest-loss”: the value from part-loading the nuclear unit in order to reduce the need for FR. Note that this value of $P_\textrm{L}$ is set by the power output of the nuclear unit in this system, since the nuclear plant drives the largest power infeed.

The largest possible outage can now be reduced by part-loading the nuclear unit: although this unit provides the cheapest energy, by part-loading it a lower amount of FR is needed, and therefore lower overall system costs in the simultaneous clearing of the energy and FR markets may be achieved.

\begin{table}[H]
\renewcommand{\arraystretch}{1.2}
\centering
\small
\begin{tabular}{p{3.3cm}| l| p{1.2cm}| p{1.2cm}}
    \multicolumn{1}{c|}{} & Nuclear & Gen. type 1 & Gen. type 2\\
\hline
Power range (MW) & \textbf{90}-100 & 0-400 & 0-300\\ 
Number of Units & 1 & 5 & 5\\ 
FR capacity (MW) & 0 & 200 & 150\\ 
Energy bid (\pounds/MWh) & 15 & 19 & 18\\ 
Inertia constant (s) & 6 & 6 & 6\\ 
FR delivery time (s) & N/A & 7 & 10\\ \hline
\end{tabular}
\caption{\label{TableThermal_ED_PLoss1}Characteristics of thermal generators. The only different parameter with respect to Table \ref{TableThermal_ED_FR2} is given in bold font.}
\end{table}

The results are presented in Table \ref{TableResults_ED_PLoss1}:

\begin{table}[H]
\renewcommand{\arraystretch}{1.2}
\centering
\small
\begin{tabular}{p{4.13cm}| l| p{0.95cm}| p{0.95cm}}
    \multicolumn{1}{c|}{} & Nuclear & Gen. type 1 & Gen. type 2\\
\hline
Power produced (MW) & 93 & 7 & 300\\ 
FR provided (MW) & 0 & 225 & 0\\ 
Revenue from energy (\pounds) & 1767 & 133 & 5700\\
Revenue from reduced-$P_\textrm{L}$ (\pounds) & 28 & - & -\\
Revenue from FR (\pounds) & 0 & 186.8 & 0\\ 
Profit (\pounds) & 400 & 186.8 & 300\\ \hline
Energy price (\pounds/MWh) & \multicolumn{3}{c}{19}\\ \hline
Reduced-$P_\textrm{L}$ price (\pounds/MW) & \multicolumn{3}{c}{4}\\ \hline
$\textrm{FR}_1$ price (\pounds/MW) & \multicolumn{3}{c}{0.83}\\ \hline
$\textrm{FR}_2$ price (\pounds/MW) & \multicolumn{3}{c}{0.58}\\ \hline
\end{tabular}
\caption{\label{TableResults_ED_PLoss1}Solution of the frequency-secured ED for the system in Table \ref{TableThermal_ED_PLoss1}, for a high-demand case of $\textrm{P}_\textrm{D}=400\textrm{MW}$.}
\end{table}

The nadir constraint enforced in this case is:
\begin{equation} \label{nadirPLossCaseStudy}
H\cdot \left(\frac{R_1}{\textrm{T}_1} + \frac{R_2}{\textrm{T}_2} \right) \geq 
\frac{(P_{\textrm{L}})^2\cdot f_0}{4 \Delta f_{\textrm{max}}}
\end{equation}

The results in Table \ref{TableResults_ED_PLoss1} demonstrate that large nuclear units would recover their opportunity cost from being part-loaded, since the price for this service obtained from the vector dual-variables of the SOC constraint (\ref{nadirPLossCaseStudy}) is \pounds4/MW, which is the difference between the energy price and the nuclear's energy bid. Generators type 1 still make a profit from providing fast FR in 7s.

A note of caution is necessary, since the payment for a reduced loss might introduce a perverse incentive for large units. To illustrate this point, let's now consider a case for which the marginal-value of a reduced-$P_\textrm{L}$ is actually higher than the opportunity cost for the nuclear unit due to being part-loaded. Consider the system in Table \ref{TableThermal_ED_PLoss2}:

\begin{table}[H]
\renewcommand{\arraystretch}{1.2}
\centering
\small
\begin{tabular}{p{3.3cm}| l| p{1.2cm}| p{1.2cm}}
    \multicolumn{1}{c|}{} & Nuclear & Gen. type 1 & Gen. type 2\\
\hline
Power range (MW) & \textbf{95}-100 & 0-400 & 0-300\\ 
Number of Units & 1 & 5 & 5\\ 
FR capacity (MW) & 0 & 200 & 150\\ 
Energy bid (\pounds/MWh) & 15 & 19 & 18\\ 
Inertia constant (s) & 6 & 6 & 6\\ 
FR delivery time (s) & N/A & 7 & 10\\ \hline
\end{tabular}
\caption{\label{TableThermal_ED_PLoss2}Characteristics of thermal generators. The only different parameter with respect to Table \ref{TableThermal_ED_PLoss1} is given in bold font.}
\end{table}

The results for this system are shown in Table \ref{TableResults_ED_PLoss2}:

\begin{table}[H]
\renewcommand{\arraystretch}{1.2}
\centering
\small
\begin{tabular}{p{4.13cm}| l| p{0.95cm}| p{0.95cm}}
    \multicolumn{1}{c|}{} & Nuclear & Gen. type 1 & Gen. type 2\\
\hline
Power produced (MW) & 95 & 19.3 & 285.7\\ 
FR provided (MW) & 0 & 225 & 14.3\\ 
Revenue from energy (\pounds) & 1805 & 366.7 & 5428.3\\
Revenue from reduced-$P_\textrm{L}$ (\pounds) & 35.5 & - & -\\
Revenue from FR (\pounds) & 0 & 315 & 14.3\\ 
Profit (\pounds) & 415.5 & 315 & 300\\ \hline
Energy price (\pounds/MWh) & \multicolumn{3}{c}{19}\\ \hline
Reduced-$P_\textrm{L}$ price (\pounds/MW) & \multicolumn{3}{c}{7.1}\\ \hline
$\textrm{FR}_1$ price (\pounds/MW) & \multicolumn{3}{c}{1.4}\\ \hline
$\textrm{FR}_2$ price (\pounds/MW) & \multicolumn{3}{c}{1}\\ \hline
\end{tabular}
\caption{\label{TableResults_ED_PLoss2}Solution of the frequency-secured ED for the system in Table \ref{TableThermal_ED_PLoss2}, for a high-demand case of $\textrm{P}_\textrm{D}=400\textrm{MW}$.}
\end{table}
All the part-loading allowed is used in this case, so the nuclear unit would make a higher profit from part-loading than if it were to sell all the energy possible. If the nuclear unit would collect the price of this ``reduced-loss" service, a perverse incentive would be created: the nuclear plant could benefit by making things worse for the whole system, i.e. by aiming for a high MSG that would increase the need for frequency services. One way of doing so would be for example to submit a higher-than-truthful MSG to the system operator.

This unwanted incentive can be avoided if the nuclear unit only gets compensated for the opportunity cost of not selling all available power in the energy market, instead of being compensated from the marginal value of $P_\textrm{L}$. The value of this opportunity cost can be easily calculated from the difference between the marginal price for energy and the marginal cost of the nuclear unit. By collecting this revenue, nuclear units do not need to receive an out-of-market payment, since they would recover their high investment just as planned (as they were conceived as must-run units providing base load at all times) but never make a perverse profit.

Note that the value of nuclear part-loading in a real system would be limited by the ramp rates of the nuclear fleet, since typically a large nuclear unit would have to be part-loaded several hours in advance, as shown in our previous work \cite{LuisEFR}. We have limited ourselves to a single-period market clearing in these examples for the sake of clarity, but the marginal-pricing formulation presented here could be directly applied to a multi-period market clearing.

\subsection{Value and incentives for inertia} \label{SectionInertia}

In this section we run several frequency-constrained UCs, so that the co-optimisation of inertia along with other frequency services can be analysed. System inertia is inherently related to the commitment state of thermal generators, and can be calculated with the following expression:
\begin{equation}
H=\sum_{g \in \mathcal{G}}\textrm{H}_g\cdot \textrm{P}_g^{\textrm{max}}\cdot y_g - \textrm{P}^{\textrm{max}}_{\textrm{L}}\cdot \textrm{H}_\textrm{L}
\end{equation}
Which considers that the largest unit won't contribute to system inertia after an outage. 

The formulation of the single-period UC with frequency constraints used here is as follows:
\begin{alignat}{3}
& {\text{min}} \quad
&& \sum_{g \in \mathcal{G}} \textrm{c}^\textrm{nl}_g\cdot y_g + \textrm{c}^\textrm{m}_g \cdot P_g  \label{ObjectiveFunction_UC}\\
& \text{s.t.} \quad
&& \sum_{g \in \mathcal{G}}P_g + \textrm{P}_\textrm{R} -P^\textrm{curt}_\textrm{R} = \textrm{P}_\textrm{D} \label{LoadBalance_UC}\\
& && y_g\cdot\textrm{P}^\textrm{msg}_g \leq P_g \leq y_g\cdot\textrm{P}^\textrm{max}_g && \quad \forall g \in \mathcal{G} \label{GenLimits_UC}\\
& && 0 \leq R_g \leq y_g\cdot\textrm{R}^\textrm{max}_g && \quad \forall g \in \mathcal{G} \label{FR_Limits_UC}\\
& && 0 \leq R_g \leq \textrm{P}^\textrm{max}_g-P_g && \quad \forall g \in \mathcal{G} \label{FR_Limits2_UC}\\
& && 0 \leq P^\textrm{curt}_\textrm{R} \leq \textrm{P}_\textrm{R}  \label{REScurtailed_limit_UC}\\
& && \textrm{RoCoF} \; \textrm{constraint} \nonumber\\
& && \textrm{Nadir} \; \textrm{constraints} \nonumber\\
& && \textrm{q-s-s} \; \textrm{constraint} \nonumber
\end{alignat}
Constraint (\ref{LoadBalance_UC}) enforces the power balance, (\ref{GenLimits_UC}) enforce generation limits, (\ref{FR_Limits_UC}) enforce FR-provision limits, (\ref{FR_Limits2_UC}) limits FR provision to be below the headroom of each generator, and constraint (\ref{REScurtailed_limit_UC}) limits RES curtailment to be below RES generation. In this UC we consider the economic inefficiencies of a part-loaded generator through their MSG and no-load cost.

As explained in Section \ref{SectionPricing}, the procedure to apply the marginal-pricing scheme and overcome the integrality of the UC is: 1) Solve the Mixed-Integer Second-Order Cone UC; 2) Solve again, relaxing the commitment decision for generators from binary to continuous and enforcing only the nadir constraint chosen in step 1, so that the problem becomes convex and therefore marginal prices can be computed. The marginal prices can be obtained for the dual variables of the problem in step 2, while the feasible operating solution is given by the original MISOCP problem in step 1.

The details of the system used for these case studies are shown in Table \ref{TableThermal_UC}:
\begin{table}[H]
\renewcommand{\arraystretch}{1.2}
\centering
\small
\begin{tabular}{p{3.3cm}| l| p{1.2cm}| p{1.2cm}}
    \multicolumn{1}{c|}{} & Nuclear & Gen. type 1 & Gen. type 2\\
\hline
Power range  per unit (MW) & 1800-1800 & 250-500 & 75-150\\ 
Number of Units & 1 & 30 & 30\\  
No-load cost $\textrm{c}^{\textrm{nl}}_g$ (\pounds) & 0 & 500 & 500\\ 
Marginal cost $\textrm{c}^{\textrm{m}}_g$ (\pounds/MWh) & 10 & 95 & 50\\ 
Inertia constant (s) & 5 & 5 & 5\\ 
FR delivery time (s) & N/A & 7 & 10\\ \hline
\end{tabular}
\caption{\label{TableThermal_UC}Characteristics of thermal generators. }
\end{table}

The solution of the frequency-secured UC for a low-RES-generation case is given in Table \ref{TableResults_UC1}. RES are assumed to have a zero marginal cost (i.e. $\textrm{c}_g^\textrm{m}$=\pounds0/\textrm{MWh}).

\begin{table}[H]
\renewcommand{\arraystretch}{1.2}
\centering
\small
\begin{tabular}{p{4.11cm}| l| p{0.96cm}| p{0.96cm}}
    \multicolumn{1}{c|}{} & Nuclear & Gen. type 1 & Gen. type 2\\
\hline
Online units & 1 & 24 & 30\\
Power produced (GW) & 1.8 & 7.7 & 4.5\\ 
FR provided (GW) & 0 & 4.3 & 0\\ 
Operating cost (\pounds k) & 18 & 743.5 & 240\\
Revenue from energy (\pounds k) & 172.4 & 737.6 & 431\\
Profit from energy (\pounds k) & 154.4 & -5.9 & 191\\ 
Total revenue (\pounds k) & 172.8 & 743.4 & 431.9\\
Total profit (\pounds k) & 154.8 & -0.1 & 191.9\\ \hline
Energy price (\pounds/MWh) & \multicolumn{3}{c}{95.79}\\ \hline
Inertia price (\pounds/MWs) & \multicolumn{3}{c}{0.041}\\ \hline
$\textrm{FR}_1$ price (\pounds/MW) & \multicolumn{3}{c}{0.79}\\ \hline
$\textrm{FR}_2$ price (\pounds/MW) & \multicolumn{3}{c}{0.55}\\ \hline
RES accommodated (GW) & \multicolumn{3}{c}{10GW (0GW curtailed)}\\ \hline
\end{tabular}
\caption{\label{TableResults_UC1}Solution of the frequency-secured UC for the system in Table \ref{TableThermal_UC}, for a demand case of $\textrm{P}_\textrm{D}=24\textrm{GW}$ and a low-RES-generation case of $\textrm{P}_\textrm{R}=10\textrm{GW}$.}
\end{table}

The enforced nadir constraint in these cases is:
\begin{equation} 
H\cdot \left(\frac{R_1}{\textrm{T}_1} + \frac{R_2}{\textrm{T}_2}\right) \geq 
\frac{(P_{\textrm{L}})^2\cdot f_0}{4 \Delta f_{\textrm{max}}}
\end{equation}

These results illustrate that the prices for frequency services would be small in a low-RES-generation condition, since the thermal units are online to cover the net demand (demand minus RES-generation) and therefore frequency services become a by-product of energy production. It is interesting to notice that a make-whole payment would be necessary for Generators type 1, and this make-whole payment would be much higher if generators were not compensated for inertia and FR (the difference in make-whole payments can be obtained from comparing the profit from energy with the total profit). Although pricing uplift, which is the cause for make-whole payments, is out of the scope of this work, this example demonstrates that remunerating generators for these ancillary services would reduce the amount of make-whole payments, a conclusion that was also reached in \cite{ElaII}.

Now let's analyse the incentives for higher inertia provision, by increasing the inertia constant of Generators type 2 to 6s:
\begin{table}[H]
\renewcommand{\arraystretch}{1.2}
\centering
\small
\begin{tabular}{p{4.11cm}| l| p{0.96cm}| p{0.96cm}}
    \multicolumn{1}{c|}{} & Nuclear & Gen. type 1 & Gen. type 2\\
\hline
Online units & 1 & 24 & 30\\
Power produced (GW) & 1.8 & 7.7 & 4.5\\ 
FR provided (GW) & 0 & 4.1 & 0\\ 
Operating cost (\pounds k) & 18 & 743.5 & 240\\
Revenue from energy (\pounds k) & 172.4 & 737.7 & 431.1\\
Profit from energy (\pounds k) & 154.4 & -5.8 & 191.1\\ 
Total revenue (\pounds k) & 172.8 & 743.4 & 432.2\\
Total profit (\pounds k) & 154.8 & -0.1 & 192.2\\ \hline
Energy price (\pounds/MWh) & \multicolumn{3}{c}{95.80}\\ \hline
Inertia price (\pounds/MWs) & \multicolumn{3}{c}{0.039}\\ \hline
$\textrm{FR}_1$ price (\pounds/MW) & \multicolumn{3}{c}{0.81}\\ \hline
$\textrm{FR}_2$ price (\pounds/MW) & \multicolumn{3}{c}{0.56}\\ \hline
RES accommodated (GW) & \multicolumn{3}{c}{10GW (0GW curtailed)}\\ \hline
\end{tabular}
\caption{\label{TableResults_UC2}Solution of the frequency-secured UC for the system in Table \ref{TableThermal_UC}, but considering that the inertia constant of Gen. type 2 is increased to 6s. The demand level is set to $\textrm{P}_\textrm{D}=24\textrm{GW}$ and a low-RES-generation case of $\textrm{P}_\textrm{R}=10\textrm{GW}$ is considered, same as for the results shown in Table \ref{TableResults_UC1}.}
\end{table}

The results show that Generators type 2 slightly increase their revenue due to increasing their inertia constant, but this increase is indeed very small, of roughly 0.1\% of the total profit. However, the incentives for inertia provision will become clear in a high-RES-generation case, for which the results are presented in Table \ref{TableResults_UC3}:

\begin{table}[H]
\renewcommand{\arraystretch}{1.2}
\centering
\small
\begin{tabular}{p{3.8cm}| l| p{0.96cm}| p{0.96cm}}
    \multicolumn{1}{c|}{} & Nuclear & Gen. type 1 & Gen. type 2\\
\hline
Online units & 1 & 17 & 27\\
Power produced (GW) & 1.8 & 4.25 & 2.05\\ 
FR provided (GW) & 0 & 4.25 & 2\\ 
Operating cost (\pounds k) & 18 & 412.2 & 116\\
Revenue from energy (\pounds k) & 0 & 0 & 0\\
Profit from energy (\pounds k) & -18 & -412.2 & -116\\ 
Total revenue (\pounds k) & 41.4 & 413.1 & 166.7\\
Total profit (\pounds k) & 23.4 & 0.9 & 50.7\\ \hline
Energy price (\pounds/MWh) & \multicolumn{3}{c}{0}\\ \hline
Inertia price (\pounds/MWs) & \multicolumn{3}{c}{4.6}\\ \hline
$\textrm{FR}_1$ price (\pounds/MW) & \multicolumn{3}{c}{51.2}\\ \hline
$\textrm{FR}_2$ price (\pounds/MW) & \multicolumn{3}{c}{35.9}\\ \hline
RES accommodated (GW) & \multicolumn{3}{c}{15.9GW (2.1GW curtailed)}\\ \hline
\end{tabular}
\caption{\label{TableResults_UC3}Solution of the frequency-secured UC for the system in Table \ref{TableThermal_UC}, for a demand case of $\textrm{P}_\textrm{D}=24\textrm{GW}$ and a high-RES-generation case of $\textrm{P}_\textrm{R}=18\textrm{GW}$.}
\end{table}

Note that in this high-RES-generation case the price of energy drops to \pounds0/MWh due to RES curtailment, since the next MWh of energy could be provided by RES for free. Therefore, all the operating cost of the system is related to the cost of frequency services. The price of inertia has increased by more than a factor 100 compared to the low-RES condition in Table \ref{TableResults_UC1}, while the prices of $\textrm{FR}_1$ and $\textrm{FR}_2$ have increased by a factor 65. 

Table \ref{TableResults_UC3} also shows that all generating units are operating at MSG, which illustrates the main driver for RES curtailment in a low-inertia system: since thermal generators providing inertia and FR must be online generating at least MSG, the sum of MSG-energy from all these sources results in RES curtailment of an equivalent value. Regarding the need for a make-whole payment for Generators type 1, it has been completely eliminated by the compensation for inertia and FR.

Let's again increase the inertia constant of Generators type 2 to 6s, to understand the incentives for higher inertia provision:
\begin{table}[H]
\renewcommand{\arraystretch}{1.2}
\centering
\small
\begin{tabular}{p{3.8cm}| l| p{0.96cm}| p{0.96cm}}
    \multicolumn{1}{c|}{} & Nuclear & Gen. type 1 & Gen. type 2\\
\hline
Online units & 1 & 16 & 28\\
Power produced (GW) & 1.8 & 4 & 2.1\\ 
FR provided (GW) & 0 & 4 & 2.1\\ 
Operating cost (\pounds k) & 18 & 388 & 119\\
Revenue from energy (\pounds k) & 0 & 0 & 0\\
Profit from energy (\pounds k) & -18 & -388 & -119\\ 
Total revenue (\pounds k) & 39.6 & 388.4 & 188.8\\
Total profit (\pounds k) & 21.6 & 0.4 & 69.8\\ \hline
Energy price (\pounds/MWh) & \multicolumn{3}{c}{0}\\ \hline
Inertia price (\pounds/MWs) & \multicolumn{3}{c}{4.4}\\ \hline
$\textrm{FR}_1$ price (\pounds/MW) & \multicolumn{3}{c}{53.1}\\ \hline
$\textrm{FR}_2$ price (\pounds/MW) & \multicolumn{3}{c}{37.1}\\ \hline
RES accommodated (GW) & \multicolumn{3}{c}{16.1GW (1.9GW curtailed)}\\ \hline
\end{tabular}
\caption{\label{TableResults_UC4}Solution of the frequency-secured UC for the system in Table \ref{TableThermal_UC}, but considering that the inertia constant of Gen. type 2 is increased to 6s. The demand level is set to $\textrm{P}_\textrm{D}=24\textrm{GW}$ and a high-RES-generation case of $\textrm{P}_\textrm{R}=18\textrm{GW}$ is considered, same as for the results shown in Table \ref{TableResults_UC3}.}
\end{table}
The profit for Generators type 2 has increased by 38\% due to increasing their inertia constant by 1s. Furthermore, RES curtailment has been reduced by 0.2GW due to this extra inertia, when compared to the results in Table \ref{TableResults_UC3}.

In conclusion, our pricing scheme promotes higher provision of inertia, which can greatly increase the revenue for generators during high-RES-generation conditions. In turn, the whole system would benefit from a lower RES curtailment, which would reduce overall costs. This also promotes investment in RES that provide services such as synthetic inertia, as they could achieve a revenue even during periods of RES curtailment when the energy price drops to \pounds0/MWh, potentially eliminating the need for feed-in tariffs since the uncertainty in revenues for RES would be reduced.

\section{Conclusions and Future Work} \label{SectionConclusions}

Following the deduction of Mixed-Integer Second-Order Cone Programming (MISOCP) frequency-security constraints for low-inertia power systems, this paper has proposed a marginal-pricing scheme for frequency services such as inertia, different types of FR and a reduced largest-loss. Through several simple case studies, we demonstrate that the proposed marginal-pricing methodology for frequency services puts in place the right incentives: providers have been shown to be incentivised to increase their inertia and achieve faster delivery of FR, as they are rewarded according to their speed of delivery. In addition, this methodology could promote investment in Renewable Energy Sources (RES) that provide frequency services, as they would collect revenue even during times of RES curtailment. We further demonstrate that large nuclear units would recover their opportunity cost from being part-loaded, but we also warn against possible perverse incentives and how to mitigate them.

Three main lines for future work can be identified. First, a methodology should be defined to map the complex Frequency Response (FR) dynamics in a real power system, driven usually by droop controllers, into a ramp defined by the three parameters $R_i$, $\textrm{T}_i$ and $\textrm{T}_{\textrm{del,}i}$. This methodology would allow providers of frequency services to simply submit these three parameters to the system operator, who then could simultaneously clear the energy and FR markets using our proposed pricing formulation. Second, integrating the contribution from load damping into the frequency-security constraints should be investigated, as previous works have identified that load damping can provide important economic savings for the provision of frequency services \cite{LuisEFR}. Third, the implications of other uplift-pricing methodologies should be investigated, as the present work has considered the dispatchable model, which relaxed the commitment decision of generators, in order to obtain marginal prices for frequency services.

\section*{}

\bibliography{Luis_PhD.bib}

\end{document}

%% file: FrequencySimulation.tikz
%
%
\definecolor{mycolor1}{rgb}{1.00000,0.26275,0.26275}%
\begin{tikzpicture}

\begin{axis}[%
axis lines = left,
width=2.65in,
height=1in,
clip=false, 
at={(0in,0in)},
scale only axis,
y label style={at={(axis description cs:-0.1,0.3)},anchor=west},
xmin=0,
xmax=16,
xlabel style={font=\footnotesize},
xlabel={time (s)},
xtick={0,5,10,15},
ytick={0,-0.772},
yticklabels={0,$\Delta f_\textrm{max}$},
ymin=-0.9,
ymax=0.1,
ylabel style={font=\footnotesize},
ylabel style={align=center},
ylabel={$\Delta f$ (Hz)},
axis background/.style={fill=white},
xticklabel style={font=\footnotesize},
yticklabel style={font=\footnotesize},
y axis line style = {-} 
]






\addplot [color=black, line width=1.0pt]
  table[row sep=crcr]{%
0	0\\
0.00120972051726858	-0.000418719516076328\\
0.00362916155180574	-0.00125609745291386\\
0.00604860258634289	-0.00209335064662660\\
0.00908408191788466	-0.00314360097412492\\
0.0222765953549888	-0.00770546147395511\\
0.0354691087920929	-0.0122626580524518\\
0.0486616222291970	-0.0168146635124410\\
0.0618541356663011	-0.0213610643981047\\
0.0788125064385063	-0.0271964896180471\\
0.0957708772107115	-0.0330215929553642\\
0.116728616427940	-0.0402055843765424\\
0.137686355645169	-0.0473726282928813\\
0.164145174046676	-0.0563960291866439\\
0.190603992448184	-0.0653910933990979\\
0.226202595688492	-0.0774477428847648\\
0.261801198928800	-0.0894512129696612\\
0.314092611333326	-0.106985060560713\\
0.382831763012408	-0.129852664584205\\
0.484436477311644	-0.163266663578856\\
0.500662131431751	-0.168559336000499\\
0.516887785551858	-0.173839755535063\\
0.528572172929525	-0.177634521705181\\
0.536903126606462	-0.180336117710059\\
0.557640963137047	-0.187045904164569\\
0.578378799667631	-0.193733483587922\\
0.599116636198216	-0.200398065848079\\
0.619854472728801	-0.207038892907320\\
0.640592309259386	-0.213655237667147\\
0.668517330753034	-0.222524669032511\\
0.696442352246683	-0.231346832937915\\
0.732538047095854	-0.242677687585673\\
0.768633741945024	-0.253923936463669\\
0.814067197766844	-0.267954690657296\\
0.859500653588664	-0.281841631629909\\
0.916080300169439	-0.298926982723753\\
0.972659946750214	-0.315773710492966\\
0.986243672226309	-0.319782068106393\\
0.999827397702403	-0.323776143320280\\
1.00795092751710	-0.326157872758370\\
1.01607445733179	-0.328534427544105\\
1.02212369782986	-0.330300764413893\\
1.04690638393902	-0.337506625052080\\
1.07168907004819	-0.344662863216207\\
1.10220749341100	-0.353405813189189\\
1.13272591677380	-0.362070579058557\\
1.17256679500647	-0.373261906850161\\
1.21240767323914	-0.384314255286318\\
1.26691475289110	-0.399204420489783\\
1.32142183254307	-0.413821947489023\\
1.39383062836068	-0.432807346996175\\
1.46623942417828	-0.451287260358085\\
1.55945496990198	-0.474314009275361\\
1.65267051562568	-0.496463300861004\\
1.76974458247863	-0.523012355114412\\
1.88681864933157	-0.548126063438347\\
2.02834830984333	-0.576545387096768\\
2.16987797035509	-0.602827711224187\\
2.31140763086685	-0.626971530408187\\
2.49147519512293	-0.654619365268584\\
2.67154275937901	-0.678875164944035\\
2.85161032363508	-0.699813269458074\\
3.07902465522130	-0.721677448704875\\
3.30643898680752	-0.738646223882679\\
3.58023163566359	-0.752983991817538\\
3.85402428451965	-0.761159331645620\\
4.12781693337572	-0.763752947933122\\
4.47957192688155	-0.759915566215414\\
4.83132692038738	-0.749182541855344\\
5.30343116258168	-0.726310659569557\\
5.77553540477598	-0.696315698640737\\
6.24763964697027	-0.661927523540426\\
6.71974388916457	-0.625545902594144\\
7.19184813135887	-0.589184330858562\\
7.66395237355316	-0.554453045456914\\
8.13605661574746	-0.522564772237778\\
8.60816085794176	-0.494357396390921\\
9.08026510013605	-0.470330284917576\\
9.55236934233035	-0.450689387640223\\
10.0244735845246	-0.435397177779720\\
10.4965778267189	-0.424223999848115\\
10.9686820689132	-0.416798056433025\\
11.4407863111075	-0.412651908395074\\
11.9128905533018	-0.411263975152162\\
12.3849947954961	-0.412094077661859\\
12.8570990376904	-0.414612553022278\\
13.3292032798847	-0.418322876546752\\
13.8013075220790	-0.422778052470378\\
14.2734117642733	-0.427591279162688\\
14.7455160064676	-0.432441564229529\\
};

\addplot [color=black, dashed]
table[row sep=crcr]{%
0 -0.772\\
15 -0.772\\
};

\addplot [color=black, dashed]
table[row sep=crcr]{%
-0.8 0.32\\
1.5 -0.6\\
};

\node[right, align=left]
at (axis cs:-0.6,0.29) {\scriptsize $\textrm{RoCoF}_\textrm{max}$};

\end{axis}
\end{tikzpicture}%

%% file: FR_services.tikz
%
%
\definecolor{mycolor1}{rgb}{1.00000,0.26275,0.26275}%
\definecolor{mycolor2}{rgb}{0.1953125,0.80078125,0.1953125}
\definecolor{mycolor3}{rgb}{0.578125,0,0.82421875}
\begin{tikzpicture}

\begin{axis}[%
axis lines = left, 
width=2.8in,
height=1.4in,
at={(0in,0in)},
scale only axis,
clip=false,
xmin=0,
xmax=3+0.15,
xtick={0,0.3,1,2},
xticklabels={{},{}},
ymin=0,
ymax=415,
ytick={0,100,125,150,375},
yticklabels={{},{}},
axis background/.style={fill=white},
title style={font=\Huge},xlabel style={font={\color{blue}\bfseries}},ylabel style={font=\tiny},legend style={font=\scriptsize},ticklabel style={font=\color{red}},
legend style={at={(0.65,0.85)}, anchor=north west, legend cell align=left, align=left, draw=white!15!black},legend style={font=\footnotesize},
legend style={fill opacity=0,text opacity=1,draw=none} 
]

\addplot [color=black, dashed, line width=1.0pt]
  table[row sep=crcr]{%
0	0\\
0.3	163.75\\
1	312.5\\
2	375\\
3	375\\
}; \addlegendentry{Total FR}

\addplot [color=black, line width=1.0pt]
  table[row sep=crcr]{%
0	0\\
0.3	100\\
3	100\\
}; \addlegendentry{$\mbox{FR}_1$}

\addplot [color=mycolor1, line width=1.0pt]
  table[row sep=crcr]{%
0	0\\
2	125\\
3	125\\
}; \addlegendentry{$\mbox{FR}_2$}

\addplot [color=mycolor2, line width=1.0pt]
  table[row sep=crcr]{%
0	0\\
1	150\\
3	150\\
}; \addlegendentry{$\mbox{FR}_3$}



\node[left, align=left] (A)
at (axis cs:0,-66) {};
\node[right, align=left] (B)
at (axis cs:0.3,-66) {};
\draw [<->,>=stealth] (A) -- (B);
\node[right, align=left]
at (axis cs:0.015,-34) {\small $\textrm{T}_1$};

\node[left, align=left] (C)
at (axis cs:0,-92) {};
\node[right, align=left] (D)
at (axis cs:1,-92) {};
\draw [<->,>=stealth] (C) -- (D);
\node[right, align=left]
at (axis cs:0.39,-62) {\small $\textrm{T}_2$};

\node[left, align=left] (C)
at (axis cs:0,-118) {};
\node[right, align=left] (D)
at (axis cs:2,-118) {};
\draw [<->,>=stealth] (C) -- (D);
\node[right, align=left]
at (axis cs:1.19,-88) {\small $\textrm{T}_3$};

\node[right, align=left]
at (axis cs:-0.35,96) {\small $R_1$};
\node[right, align=left]
at (axis cs:-0.35,128) {\small $R_3$};
\node[right, align=left]
at (axis cs:-0.35,160) {\small $R_2$};
\node[right, align=left]
at (axis cs:-0.55,375) {\small $\sum\limits_{i=1}^3 R_i$};

\node[right, align=left, rotate=90]
at (axis cs:-0.55,100) {\small FR $\; (\mbox{MW})$};
\node[right, align=left]
at (axis cs:2.65,-45) {\small $t \; (\mbox{s})$};

\end{axis}
\end{tikzpicture}%

%% file: FR_services2.tikz
%
%
\definecolor{mycolor1}{rgb}{1.00000,0.26275,0.26275}%
\definecolor{mycolor2}{rgb}{0.1953125,0.80078125,0.1953125}
\definecolor{mycolor3}{rgb}{0.578125,0,0.82421875}
\begin{tikzpicture}

\begin{axis}[%
axis lines = left, 
width=2.8in,
height=1.4in,
at={(0in,0in)},
scale only axis,
clip=false,
xmin=0,
xmax=3+0.15,
xtick={0,0.5,1,1.7,2.3},
xticklabels={{},{}},
ymin=0,
ymax=280,
ytick={0,100,150,250},
yticklabels={{},{}},
axis background/.style={fill=white},
title style={font=\Huge},xlabel style={font={\color{blue}\bfseries}},ylabel style={font=\tiny},legend style={font=\scriptsize},ticklabel style={font=\color{red}},
legend style={at={(0.065,0.9)}, anchor=north west, legend cell align=left, align=left, draw=white!15!black},legend style={font=\footnotesize},
legend style={fill opacity=0,text opacity=1,draw=none} 
]

\addplot [color=black, dashed, line width=1.0pt]
  table[row sep=crcr]{%
0.5	0\\
1	100\\
1.7	100\\
2.3	250\\
3	250\\
}; \addlegendentry{Total FR}

\addplot [color=black, line width=1.0pt]
  table[row sep=crcr]{%
0.5	0\\
1	100\\
3	100\\
}; \addlegendentry{$\mbox{FR}_1$}

\addplot [color=mycolor1, line width=1.0pt]
  table[row sep=crcr]{%
1.7	0\\
2.3	150\\
3	150\\
}; \addlegendentry{$\mbox{FR}_2$}



\node[left, align=left] (A)
at (axis cs:0,-46) {};
\node[right, align=left] (B)
at (axis cs:0.5,-46) {};
\draw [<->,>=stealth] (A) -- (B);
\node[right, align=left]
at (axis cs:0.015,-24) {\small $\textrm{T}_\textrm{del,1}$};

\node[left, align=left] (A)
at (axis cs:0.5,-46) {};
\node[right, align=left] (B)
at (axis cs:1,-46) {};
\draw [<->,>=stealth] (A) -- (B);
\node[right, align=left]
at (axis cs:0.6,-24) {\small $\textrm{T}_1$};

\node[left, align=left] (C)
at (axis cs:0,-62) {};
\node[right, align=left] (D)
at (axis cs:1.7,-62) {};
\draw [<->,>=stealth] (C) -- (D);
\node[right, align=left]
at (axis cs:1.10,-40) {\small $\textrm{T}_\textrm{del,2}$};

\node[left, align=left] (C)
at (axis cs:1.7,-62) {};
\node[right, align=left] (D)
at (axis cs:2.3,-62) {};
\draw [<->,>=stealth] (C) -- (D);
\node[right, align=left]
at (axis cs:1.85,-40) {\small $\textrm{T}_2$};

\node[right, align=left]
at (axis cs:-0.35,96) {\small $R_1$};
\node[right, align=left]
at (axis cs:-0.35,160) {\small $R_2$};
\node[right, align=left]
at (axis cs:-0.55,250) {\small $R_1\hspace{-1mm}+\hspace{-1mm}R_2$};

\node[right, align=left, rotate=90]
at (axis cs:-0.55,60) {\small FR $\; (\mbox{MW})$};
\node[right, align=left]
at (axis cs:2.65,-35) {\small $t \; (\mbox{s})$};

\end{axis}
\end{tikzpicture}%